\documentclass[12pt]{article}
\usepackage{amsmath, amssymb, amsthm, enumerate}

\newtheorem{theorem}{Theorem}
\newtheorem{lemma}[theorem]{Lemma}
\newtheorem{proposition}[theorem]{Proposition}
\newtheorem{remark}[theorem]{Remark}

\newtheorem{cor}[theorem]{Corollary}

\def\rr{{\mathbb R}}

\def\am{^{-1}}

\def\cc{{\mathbb C}}

\def\zz{{\mathbb Z}}

\def\qq{{\mathbb Q}}

\def\dd{{\cal D}}

\def\zk{{\zz}^k}

\def\om{\omega}
\def\rn{{\rr}^n}

\def\rp{{\rr}^p}
\def\su{\subset}

\def\se{\setminus}

\def\al{\alpha}

\def\ga{\gamma}
\def\Ga{\Gamma}
\def\de{\delta}
\def\De{\Delta}
\def\ep{\varepsilon}
\def\va{\vartheta}

\def\Ph{\Phi}

\def\la{\lambda}
\def\La{\Lambda}

\def\om{\omega}

\def\mmm{{\cal M}}
\def\pp{{\cal P}}
\def\cd{\cdot}
\def\stb{,\ldots ,}

\def\emp{\emptyset}

\def\V{\Vert}
\def\msk{\medskip}
\def\bsk{\bigskip}
\def\noi{\noindent}
\def\akkor{\Longrightarrow}

\def\ol{\overline}

\def\sumin{\sum_{i=1}^n}

\def\sumis{\sum_{i=1}^s}

\def\sumiN{\sum_{i=1}^N}

\def\sumjm{\sum_{j=1}^m}
\def\sumj0m{\sum_{j=0}^m}
\def\sumik{\sum_{i=1}^k}

\def\sumi0n{\sum_{i=0}^n}

\begin{tiny}\end{tiny}

\def\proof{\noi {\bf Proof.} }

\def\x1n{x_1 \stb x_n}
\def\y1n{y_1 \stb y_n}

\def\inter{{\rm int}\, }
\def\cl{{\rm cl}\, }
\def\deg{{\rm deg}\, }
\def\det{{\rm det}\, }
\def\dim{{\rm dim}\, }

\begin{document}

\title{Vector valued polynomials, exponential polynomials and vector valued
harmonic analysis}

\author{M. Laczkovich}

\footnotetext[2]{{\bf MR subject classification:} primary 39B52, secondary
22A20}
\footnotetext[3]{Supported by the Hungarian National
Foundation for Scientific Research, Grant No. K124749.}

\maketitle

\begin{abstract}
Let $G$ be a topological Abelian semigroup with unit, let $E$ be a Banach
space, and let $C(G,E)$ denote the set of continuous functions $f\colon G\to E$.
A function $f\in C(G,E)$ is a generalized
polynomial, if there is an $n\ge 0$ such that $\De _{h_1} \ldots \De _{h_{n+1}} f=0$ for every $h_1 \stb h_{n+1} \in G$, where $\De _h$ is the difference operator.
We say that $f\in C(G,E)$ is a polynomial, if it is a
generalized polynomial, and the linear span of its translates is of finite
dimension; $f$ is a w-polynomial, if $u\circ f$ is a polynomial for every $u\in
E^*$, and $f$ is a local polynomial, if it is a polynomial on every
finitely generated subsemigroup. We show that each of the classes of
polynomials, w-polynomials, generalized polynomials, local polynomials is
contained in the next class. If $G$ is an Abelian group and has a 
dense subgroup with finite torsion free rank, then these classes coincide.

We introduce the classes of exponential polynomials and
w-expo\-nential polynomials as well, establish their representations and
connection with polynomials and w-polynomials.

We also investigate spectral synthesis and analysis in the class $C(G,E)$.
It is known that if $G$ is a compact Abelian group and $E$ is a Banach space,
then spectral synthesis holds in $C(G,E)$. On the other hand, we show that if
$G$ is an infinite and discrete Abelian group and
$E$ is a Banach space of infinite dimension, then even spectral analysis
fails in $C(G,E)$. If, however, $G$ is discrete, has finite torsion free rank
and if $E$ is a Banach space of finite dimension, then spectral synthesis holds
in $C(G,E)$.
\end{abstract}

\section{Vector valued polynomials and exponential polynomials}

Let $G$ be a topological Abelian semigroup with unit. We denote the semigroup
operation by addition, and denote the unit by $0$. Let $E$ be a Banach
space over the complex field $\cc$. We denote by $E^G$ the set of maps from $G$
into $E$. 
Let $T_g$ denote the translation operator on $E^G$ defined by $T_g f(x)=f(x+g)$
for every $f\in E^G$ and $g,x\in G$. A subset $V\su E^G$ is {\it translation
invariant}, if $T_g f\in V$ whenever $f\in V$ and $g\in G$. If $f\in E^G$, then
$L_f$ denotes the linear span of $\{ T_g f \colon g\in G\}$.


The operator $\De _g$ is defined by $\De _g =T_g -T_0$. That is,
we have $\De _g f(x)=f(x+g)-f(x)$ for every $f\in E^G$ and $x\in G$.

We say that a continuous function $f\in E^G$ is a {\it generalized
polynomial},\footnote{Our terminology differs from that of \cite{Sz2}.}
if there is an $n\ge 0$ such that $\De _{h_1} \ldots \De _{h_{n+1}} f=0$ for every
$h_1 \stb h_{n+1} \in G$. The smallest $n$ with this property is the {\it degree
of $f$}, denoted by $\deg f$. The degree of the identically zero function is
$-1$.

By Djokovi\'c's theorem \cite{Dj} (see also \cite[Section 6]{108}), a continuous
function $f\in E^G$ is a
generalized polynomial if and only if $f=\sumin f_i$, where $f_i$ is a monomial
of degree $i$ for every $i=1\stb n$, and $f_0$ is constant. By a monomial
of degree $i$ we mean a function of the form $A(x\stb x)$, where $A(x_1 \stb
x_i )$ is a map from $G^i$ to $E$ which is symmetric, and $i$-additive; that
is, additive in each variable. It is easy to see that the representation
$f=\sumin f_i$ is unique. 

It is clear that the set of generalized polynomials forms a linear
subspace of $E^G$ over $\cc$.

\begin{theorem} \label{t1}
\begin{enumerate}[{\rm (i)}]
\item A continuous function $f\in E^G$ is a generalized polynomial if and
only if $u\circ f$ is a (complex valued) generalized polynomial for every
$u\in E^*$. 
\item If $f$ is a generalized polynomial, then $\deg (u\circ f)\le \deg f$
for every $u\in E^*$.
\item If $f$ is a generalized polynomial, then there is an $u\in E^*$
such that $\deg (u\circ f)= \deg f$.
\end{enumerate}
\end{theorem}
\proof
The ``only if'' direction of (i) is obvious, and so is (ii). To prove the ``if''
statement of (i), let $E_n^*$ denote the set of linear functionals $u\in E^*$
such that $\De _{h_1} \ldots \De _{h_{n+1}} (u\circ f)=0$ for every $h_1 \stb h_{n+1}
\in G$. It is easy to see that $E_n^*$ is a closed linear subspace of $E^*$ for
every $n=0,1,\ldots$.

If $u\circ f$ is a generalized polynomial for every $u\in E^*$, then $E^* =
\bigcup_{n=0}^\infty E_n^*$. Then, by the Baire category theorem, there is an $n$
such that $\inter E_n^* \ne \emp$, and thus $E_n^* =E^*$. Let $n$ be the
smallest such $n$.

We show that $\De _{h_1} \ldots \De _{h_{n+1}} f=0$ for every $h_1 \stb h_{n+1}$.
Indeed, if $\De _{h_1} \ldots \De _{h_{n+1}} f(x) \ne 0$ for some $h_1 \stb
h_{n+1}, x\in G$, then there is an $u\in E^*$ such that
$$\De _{h_1} \ldots \De _{h_{n+1}} (u\circ f)(x) =
u\left( \De _{h_1} \ldots \De _{h_{n+1}} f(x) \right) \ne 0,$$
which is impossible. This proves both (i) and (iii). \hfill $\square$

\begin{remark} \label{r1}
{\rm The continuity of the function $f$ cannot be omitted from the conditions of
the theorem. In other words, a function $f\in E^G$ such that $u\circ f$ is a
complex valued generalized polynomial, hence continuous for every
$u\in E^*$ is not necessarily continuous itself. 

As an example, let $E$ be an
infinite dimensional Banach space, and let $G$ denote the vector space of $E$
endowed with the weak topology of $E$. Then $G$ is a topological vector space.
Let $f$ denote the identity on $G$ as a map from $G$ to the Banach space $E$.
Then $f$ is not continuous, as the original topology of $E$ is strictly stronger
than the weak topology. On the other hand, if $u\in E^*$, then $u\circ f=u$
is a continuous additive function, therefore, a generalized polynomial.}
\end{remark}

We show that if $G$ is a normed linear space, then the continuity
of $f$ is a consequence of the other condition.
\begin{theorem} \label{t2}
Let $G$ be a normed linear space, and let $E$ be a Banach space. A
function $f\in E^G$ is a generalized polynomial if and
only if $u\circ f$ is a (complex valued) generalized polynomial for every
$u\in E^*$.
\end{theorem}
\proof
By Theorem \ref{t1}, we only have to show that if $u\circ f$ is a
complex valued generalized polynomial for every $u\in E^*$, then $f$ is
continuous. 

Let $f=\sum_{i=0}^n f_i$, where $f_i$ is a monomial of degree $i$ for $1\le i
\le n$, and $f_0$ is a constant. It is enough to show that $f_i$ is continuous
for every $i=1\stb n$.

It is easy to see that $f(kx)=\sum_{i=0}^n k^i \cd f_i (x)$ for every $x\in G$
and for every positive integer $k$.
These equations for $k=1\stb n+1$ constitute a linear system of equations with
unknowns $f_i (x)$ $(i=0\stb n)$. Since the determinant of this system
is nonzero (being a Vandermonde determinant), it follows that each $f_i (x)$
is a linear combination of $f(x)\stb f((n+1)x)$ with rational coefficients.

If $u\in E^*$, then $u\circ f$ is a generalized polynomial, hence continuous.
Then each of the functions $x\mapsto u(f(kx))$ $(k=1\stb n+1)$ is continuous,
and thus $u\circ f_i$, being a linear combinations of these functions, is also
continuous for every $i=1\stb n$.

Let $B_r$ denote the open ball $\{ x\in G\colon \V x \V _G <r\}$ (recall that
$G$ is a normed linear space by assumption). Let $1\le i\le n$ be fixed.
We show that $f_i (B_1 )$ is weakly bounded in $E$. Indeed, if $u\in E^*$, then
the continuity of $u\circ f_i$ implies that for a suitable positive integer $k$,
$|u(f_i (x))| <1$ for every $x\in B_{1/k}$. Therefore, if $x\in B_1$, then
$x/k\in B_{1/k}$, $|u(f_i (x/k))| <1$ and $|u(f_i (x))| <k^i$, showing that $u$
is bounded on $f_i (B_1 )$. This proves that $f_i (B_1 )$ is weakly bounded in
$E$. Since, in a Banach space, every weakly bounded set is originally bounded
\cite[3.18 Theorem]{Ru}, it follows that $\V y\V _E <K$ for every $y\in f_i
(B_1)$ with a suitable positive integer $K$.

If $\ep >0$ is given, then there is an integer $m$ such that $m>K/\ep$.
If $x\in B_{1/m}$, then
$$ \V f_i (x)\V _E  =\V m^{-i} f_i (mx) \V _E <K/m^i <\ep ,$$
proving that $f_i$ is continuous at zero. 
Now, it is known that if a monomial is continuous at one point, then it is
continuous everywhere. See \cite[Theorem 3.6]{Sz2}. Note that monomials are
``algebraic polynomials'' in the terminology of \cite{Sz2}, and that
the conditions of \cite[Theorem 3.6]{Sz2} are satisfied if $G$ is a normed
linear space and $E$ is a Banach space. Thus $f_i$ is continuous on $G$
for every $i=1\stb n$, and this is what we wanted to show. \hfill $\square$

\bsk
A function $f\in \cc ^G$ is said to be a {\it polynomial}, if there are
continuous additive functions $a_1 \stb a_n \colon G\to \cc$ and there is a
$P\in \cc [x_1 \stb x_n ]$ such that $f=P(a_1 \stb a_n )$. It is well-known that
every complex valued polynomial is a generalized polynomial. 

A continuous function $m\colon G\to \cc$ is called an {\it
exponential function}, if $m\ne 0$ and $m(x+y)=m(x)\cd m(y)$ for every $x,y\in
G$. 

A function $f\in \cc ^G$ is an {\it exponential polynomial}, if there are
polynomials $p_1 \stb p_n \colon G\to \cc$ and exponentials $m_1 \stb m_n$
such that $f=\sumin p_i \cd m_i$. It is well-known that a continuous function
$f\in \cc ^G$ is an exponential polynomial if and only if $\dim L_f <\infty$. 
(For Abelian semigroups see \cite{Mc}, for Abelian groups see \cite{Sz1}
and \cite[Theorem 10.2]{Sz2}. See also \cite{140} for the history of the
theorem and for a simple proof.)
In possession of this result the following definition seems reasonable.
We say that a continuous function $f\in E^G$ is an {\it exponential polynomial},
if $\dim L_f <\infty$. 

\begin{theorem} \label{t3}
A function $f\in E^G$ is an exponential polynomial if and only if there are
finitely many complex valued exponential polynomials $f_1 \stb f_k \in \cc ^G$
and elements $e_1 \stb e_k \in E$ such that $f=f_1 \cd e_1 +\ldots + f_k \cd
e_k$.
\end{theorem}
\proof
The ``if'' statement is clear: if $f_1 \stb f_k$ are exponential polynomials
and $f=f_1 \cd e_1 +\ldots + f_k \cd e_k$,
then $L_{f_1 \cd e_1} \stb L_{f_k \cd e_k}$ are of finite dimension, and then so is
$L_f$.

To prove the converse, suppose that $\dim L_f <\infty$. First we show that
the linear subspace $F$ of $E$ generated by $R(f)$, the range of $f$, is of
finite dimension. Suppose not. Then there are elements $x_1 \stb x_n \in G$ such
that $n>\dim L_f$, and $f(x_1 )\stb f(x_n )$ are linearly independent over
$\cc$. Now $n>\dim L_f$ implies that $T_{x_1} f \stb T_{x_n} f$ are linearly
dependent, and thus $\sumin c_i T_{x_i} f =0$ for some complex numbers $c_1
\stb c_n$, not all zero. Then $\sumin c_i f(x+x_i ) =0$ for every $x$. In
particular, putting $x=0$ we get $\sumin c_i f(x_i ) =0$, which contradicts the
fact that  $f(x_1 )\stb f(x_n )$ are linearly independent.

This proves that $\dim F<\infty$. Let $e_1 \stb e_k$ be a basis of $F$. Then
there are functions $f_1 \stb f_k \colon G\to \cc$ such that $f=f_1 \cd e_1 +
\ldots + f_k \cd e_k$ on $G$. Since $e_1 \stb e_k$ are linearly independent,
there are linear functionals $u_1 \stb u_k \in E^*$ such that $u_i (e_i)=1$
and $u_i (e_j )=0$ for every $1\le i,j\le k$, $i\ne j$. Thus $f_i =u_i \circ f$
for every $i=1\stb k$. Since $f$ is continuous, we can see that so are $f_1
\stb f_k$. Also, we have $\dim L_{f_i} =\dim L_{u_i \circ f} \le \dim L_{f} <\infty$,
and thus $f_i$ is an exponential polynomial for every
$i=1\stb k$. \hfill $\square$

The definition of polynomials cannot be generalized to the vector-valued
case, as Banach spaces are not algebras. However, the following observation
makes it clear what a reasonable generalization could be.

\begin{proposition}\label{p1}
A complex valued function $f\colon G\to \cc$ is a polynomial if and only
if $f$ is a generalized polynomial, and $\dim L_f <\infty$.
\end{proposition}
\proof
We noted already that every polynomial is a generalized polynomial. In fact,
one can prove by induction on $\deg P$ that
if $f=P(a_1 \stb a_n )$, where $P\in \cc [x_1 \stb x_n ]$ and 
$a_1 \stb a_n \colon G\to \cc$ are continuous additive functions, then
$f$ is a generalized polynomial of degree at most $\deg P$. Similarly,
$\dim L_f <\infty$ can also be proved by induction on $\deg P$. Or, we can argue
that if $f$ is a polynomial, then it is an exponential polynomial, since $f=
f\cd 1$ and the identically $1$ function is an exponential. Thus $\dim L_f <
\infty$ follows from a theorem quoted before. This proves the ``only if'' part
of the proposition.

Now suppose that $f$ is a generalized polynomial and $\dim L_f <\infty$.
The latter condition implies that $f$ is an exponential polynomial; that is,
$f=\sumik p_i \cd m_i$, where $p_1 \stb p_k$ are polynomials and $m_1 \stb m_k$
are exponentials. We may assume that $p_1 \stb p_k$ are nonzero and $m_1 \stb
m_k$ are distinct.

It is known that the representation of a function $f\colon G\to \cc$ in the form
$\sumis p_i \cd m_i$, where $m_1 \stb m_s$ are distinct exponentials and $p_1
\stb p_s$ are nonzero generalized polynomials is unique (if exists).
For Abelian groups this is proved in \cite[Lemma 4.3, p. 41]{Sz2}
and in \cite[Lemma 6]{L}. It is easy to check that the proof of
\cite[Lemma 6]{L} works in Abelian semigroups as well. The uniqueness follows
also from Lemma \ref{l2} below.

In our case $f\cd 1=\sumis p_i \cd m_i$ and thus the uniqueness of the
representation implies $s=1$, $m_1 =1$ and $f=p_1$. Thus $f$ is a polynomial,
which proves the ``if'' part of the proposition. \hfill $\square$

The proposition above motivates the following definition: a 
function $f\in E^G$ is a {\it polynomial}, if $f$ is a generalized polynomial,
and $\dim L_f <\infty$.

\begin{theorem} \label{t4}
A function $f\in E^G$ is a polynomial if and only if there are
finitely many complex valued polynomials $f_1 \stb f_k$ and elements
$e_1 \stb e_k \in E$ such that $f=f_1 \cd e_1 +\ldots + f_k \cd e_k$.
\end{theorem}
\proof
Suppose that $f=f_1 \cd e_1 +\ldots + f_k \cd e_k$, where $f_1 \stb f_k$ are
complex valued polynomials and $e_1 \stb e_k \in E$. Then $f_1 \stb f_k$ are
generalized polynomials, hence so are $f_1 \cd e_1 \stb f_k \cd e_k$, and then
so is $f$. Also, $L_{f_1} \stb L_{f_k}$ are of finite dimension, implying
$\dim L_f < \infty$. This proves the ``if'' statement.

If $f$ is a polynomial, then $\dim L_f <\infty$. By Theorem \ref{t3},
this implies that $f=f_1 \cd e_1 +\ldots + f_k \cd e_k$, where $f_1 \stb f_k$ are
complex valued exponential polynomials, and $e_1 \stb e_k \in E$. We may assume
that $e_1 \stb e_k$ are linearly independent. Then there are linear
functionals $u_1 \stb u_k \in E^*$ such that $u_i (e_i)=1$ and $u_i (e_j )=0$
for every $1\le i,j\le k$, $i\ne j$. Thus $f_i =u_i \circ f$ for every
$i=1\stb k$. Since $f$ is a generalized polynomial, we can see that so are $f_1
\stb f_k$. Summing up: $f_1 \stb f_k$ are generalized polynomials and
exponential polynomials. Therefore, they are polynomials by Proposition
\ref{p1}. \hfill $\square$

\begin{theorem} \label{t5}
A function $f\in E^G$ is an exponential polynomial if and only if
$f=\sumik m_i \cd p_i$, where $m_1 \stb m_k$ are (complex valued) exponentials,
and $p_1 \stb p_k \in E^G$ are polynomials.
\end{theorem}
\proof
If $m$ is an exponential and $p\in E^G$, then $T_g (m\cd p)= m(g)\cd m\cd
T_g p$, and thus
$$T_g (m\cd p)\in \{ m\cd \phi \colon \phi \in L_p \}$$
for every $g\in G$. If $p$ is a polynomial, then $L_p$ is of finite dimension,
and then so is $L_{m\cd p}$. Thus $m\cd p$ is an exponential polynomial whenever
$m$ is a complex valued exponential and $p$ is a polynomial. From this
observation the ``if'' part of the statement of the theorem is obvious.

To prove the ``only if'' part, let $f$ be an exponential polynomial.
By Theorem \ref{t3}, $f=f_1 \cd e_1 +\ldots + f_k \cd e_k$, where
$f_1 \stb f_k$ are complex valued exponential polynomials and
$e_1 \stb e_k \in E$. Let $f_i =\sum_{j=1}^{n_i} p_{ij}\cd m_{ij}$, where
$p_{ij}$ is a complex valued polynomial and $m_{ij}$ is an exponential for every
$j=1\stb n_i$. Then
$$f=\sumik \sum_{j=1}^{n_i} m_{ij}\cd (p_{ij} \cd e_i ),$$
where $p_{ij} \cd e_i \in E^G$ is a polynomial for every $i,j$. \hfill $\square$

\section{The classes of w-polynomials and w-expo\-nential polynomials}

In the complex valued case the continuous additive functions are automatically
polynomials. In the vector valued setting this is not the case, as the following
example shows.

Let $E$ be an infinite dimensional Banach space, and let $G$ be its additive
group with the same topology. Let $f\colon G\to E$ be the identity map on $G$.
Then $f$ is a continuous
additive function, but not a polynomial. Indeed, $L_f$ equals the set of
functions $x\mapsto cx+e$, where $c\in \cc$ and $e\in E$. In particular, $L_f$
contains the constant functions, and thus $\dim L_f = \infty$. Consequently,
$f$ is not a polynomial. On the other hand, it is clear that $u\circ f =u$ is a
(complex valued) polynomial for every $u\in E^*$. This motivates the following
definition.

Let $G$ be a topological Abelian semigroup with unit, and let $E$ be a Banach
space over the complex field $\cc$. 
We say that a continuous function $f\in E^G$ is a {\it w-polynomial}, if
$u\circ f$ is a (complex valued) polynomial for every $u\in E^*$.

We introduce one more variation on the theme of polynomials.
A continuous function $f\in E^G$ is called a {\it local polynomial}, if
the restriction of $f$ to any finitely generated subsemigroup of $G$ is a
polynomial.
\begin{theorem} \label{t6}
Consider the following properties that a continuous function $f\in E^G$ may
have:
\begin{enumerate}[{\rm (i)}]
\item $f$ is a polynomial, 
\item $f$ is a w-polynomial,
\item $f$ is a generalized polynomial,
\item $f$ is a local polynomial.
\end{enumerate}
Then we have {\rm (i)}$\akkor${\rm (ii)}$\akkor${\rm (iii)}$\akkor${\rm (iv)}.
\end{theorem}
\proof
(i)$\akkor$(ii): If $f$ is a polynomial, then $f$ is a generalized polynomial,
and $\dim L_f <\infty$. It is clear that if $u\in E^*$, then $u\circ f$ has
the same properties, and thus $u\circ f$ is a polynomial by Proposition
\ref{p1}.

\noi
(ii)$\akkor$(iii): If $f$ is a w-polynomial, then $u\circ f$ is a 
polynomial for every $u\in E^*$. Thus $u\circ f$ is a 
generalized polynomial for every $u\in E^*$. Therefore, by Theorem \ref{t1}, 
$f$ is a generalized polynomial.

\noi
(iii)$\akkor$(iv): Let $f$ be a generalized polynomial. As we mentioned earlier,
this implies, by Djokovi\'c's theorem \cite{Dj}, that $f=\sumin f_i$, where
$f_i$ is a monomial of degree $i$ for every $i=1\stb n$, and $f_0$ is constant.
Let $f_i (x)= A_i (x\stb x)$, where $A_i (x_1 \stb x_i )$ is symmetric and
additive in each variable.

Let $H$ be a finitely generated subsemigroup of $G$, and let $h_1 \stb h_k$ be a
generating system of $H$. It is clear that the restriction $f|_H$ is also a
generalized polynomial. We prove that $\dim L_{f|H} <\infty$.

If $x,y\in H$, then $f(x+y)=\sumin f_i (x+y)$. Since $A_i$ is symmetric and
additive in each variable, we have $f_i (x+y)=A_i (x+y \stb x+y)= \sum_{j=0}^i
g_i (x,y)$, where
$$g_i (x,y)={i \choose j} A(\underbrace{x \stb x}_{j} \underbrace{y \stb y}_{i-j} ) .$$
Since $y$ is a linear combinations with nonnegative integer coefficients
of the elements $h_1 \stb h_k$, it follows that $g_i (x,y)$
is a linear combinations with nonnegative integer coefficients of the functions
$A_i (x\stb x, h_{\nu _1} \stb h_{\nu _{j-i}} )$, where $\nu _1 \stb \nu _{j-i} \in
\{ 1 \stb  k\}$. Thus $T_y f_i$ is the linear combination of finitely
many functions that are independent of the choice of $y\in H$. Therefore,
$\dim L_{f_i |H} <\infty$ for every $i=1\stb n$, and thus $\dim L_{f|H} <\infty$.
This proves that $f$ is a polynomial on $H$. Since $H$ was an arbitrary
finitely generated subsemigroup of $G$, it follows that $f$ is a local
polynomial on $G$. \hfill $\square$

\bsk
Note that if $E$ is finite dimensional, then every $w$-polynomial is a
polynomial. Indeed, if $f$ is a $w$-polynomial and $e_1 \stb e_n$ is a basis of
$E$, then $f=f_1 \cd e_1 +\ldots + f_n \cd e_n$, where
$f_1 \stb f_n \in \cc ^G$. The argument of the proof of Theorem \ref{t4}
gives that $f_i =u_i \circ f$ with suitable $u_1 \stb u_n \in E^*$,
and thus $f_1 \stb f_n$ are polynomials. Thus $f$ itself is a polynomial
by Theorem \ref{t4}.

If $G$ is a finitely generated semigroup, then every local polynomial is a
polynomial,
and thus properties (i)-(iv) are equivalent. We show that for Abelian groups
somewhat more is true.
If $G$ is an Abelian group, then we denote by $r_0 (G)$ the torsion free rank
of $G$; that is, the cardinality of a maximal independent system of elements of
infinite order.

\begin{theorem} \label{t7}
Let $G$ be a topological Abelian group, and suppose that there is a dense
subgroup $H$ of $G$ such that $r_0 (H) <\infty$. Then properties
{\rm (i)-(iv)} listed in Theorem \ref{t6} are equivalent.
\end{theorem}




\proof
First we assume that $r_0 (G)<\infty$. Then there is a finitely generated
subgroup $H$ of $G$ such that the factor group $G/H$ is torsion. In other words,
for every $h\in G$ there is a positive integer $k$ such that $kh\in H$.

It is enough to prove that if $f\in E^G$ is a local polynomial, then 
$f$ is a polynomial.

First we show that $f$ is a generalized polynomial. 
Since $f$ is a local polynomial, the restriction $f|_H$ to the finitely
generated subgroup $H$ is a
polynomial. In particular, $f|_H$ is a generalized polynomial, and thus
$f|_H =\sumin f_i$, where $f_i$ is a monomial of degree $i$ for every $i=1\stb
n$, and $f_0$ is constant.

We show that $f$ is a generalized polynomial of degree $n$ on $G$. Let
$a_1 \stb a_{n+1},x \in G$ be arbitrary; we prove $\De _{a_1} \ldots
\De _{a_{n+1}} f(x)=0$.
Let $\ol H$ denote the subgroup of $G$ generated by $H$ and the elements 
$a_1 \stb a_{n+1},x$. Then $\ol H$ is finitely generated. Since $f$ is a local
polynomial, it follows that $f|_{\ol H}$ is a polynomial, hence a generalized
polynomial. Thus $f|_{\ol H} = \sumjm g_j$, where $g_i$ is a monomial of degree
$j$ for every $j=1\stb m$, and $g_0$ is constant. We may assume that $g_m$ is
not identically zero on $G$. We prove $m\le n$.

Since $H\su \ol H$, we have $f|_H =\sumjm g_j |_H$. Now the representation of
$f|_H$ as a sum of monomials is unique. If $m>n$, then necessarily
$g_m |_H =0$. Let $g_m (x)=B(x\stb x)$, where $B$ is $m$-additive.
If $h\in \ol H$, then $kh\in H$ with a suitable positive integer $k$. Then
$$g_m (h)=B(h\stb h)=k^{-m} B(kh \stb kh )= k^{-m} g_m (kh)=0,$$
as $kh \in H$. Thus $g_m$ is identically zero on $\ol H$, which is a
contradiction.

This proves $m\le n$. Then $f|_{\ol H} = \sumjm g_j$ implies that $f|_{\ol H}$
is a generalized polynomial of degree at most $n$. Since $a_1 \stb a_{n+1} ,x
\in \ol H$, we obtain $\De _{a_1} \ldots \De _{a_{n+1}} f(x)=0$, and this is what
we wanted to show. Thus $f$ is a generalized polynomial on $G$. Let 
$f|_H =\sumin f_i$, where $f_i$ is a monomial of degree $i$ for every $i=1\stb
n$, and $f_0$ is constant. 

Now we prove that $f$ is a polynomial. We only have to show that $\dim L_f <
\infty$.

Clearly, it is enough to show that $\dim L_{f_i} <\infty$ for every
$i=1\stb n$. Let $f_i (x)=A_i (x\stb x)$, where $A_i$ is symmetric and
$i$-additive.  Let $h_1 \stb h_N$ be a generating system of $H$. We prove that
$L_{f_i}$ is contained by the linear hull of the functions
\begin{equation}\label{e7}
A(\underbrace{x \stb x}_{i-j} ,h_{\nu _1}\stb  h_{\nu _j} ),
\end{equation}
where $0\le j\le i$ and $\nu _1 \stb \nu _j \in \{ 1\stb N\}$. Since the number
of these functions is finite, this will prove $\dim L_{f_i} <\infty$.

Let $h\in G$ be arbitrary. Then $T_h f_i (x)=f_i (x+h)=A_i (x+h \stb x+h)$ is a
linear combination with integer coefficients of the functions
$A(\underbrace{x \stb x}_{i-j} , h\stb h)$ $(j=0\stb i)$. Let $k$ be a positive
integer with $kh\in H$. If $j$ is fixed, then
\begin{equation}\label{e8}
A(\underbrace{x \stb x}_{i-j} ,h\stb h)=k^{-j} A(\underbrace{x \stb x}_{i-j} ,
kh\stb kh).
\end{equation}
Now $kh$, being an element of $H$, is a linear combination with integer
coefficients of the elements $h_1 \stb h_N$. Since $A_i$ is additive in each
variable, it follows from \eqref{e8} that 
$A(\underbrace{x \stb x}_{i-j} ,h\stb h)$ is a linear combination with rational
coefficients of the function listed in \eqref{e7}. This completes the proof
of $\dim L_{f_i} <\infty$.
We proved that the statement of the theorem is true if $r_0 (G)<\infty$.

Now we assume that $G$ has a dense subgroup $H$ such that $r_0 (H) <\infty$.
Suppose that $f\in E^G$ is a local polynomial. We have to prove that 
$f$ is a polynomial; that is, a generalized polynomial satisfying
$\dim L_f <\infty$.

Since $r_0 (H)<\infty$ and $f$ is a local polynomial on $H$, it follows that
$f|_H$ is a polynomial on $H$. Then $f|_H$ is a generalized
polynomial on $H$; let $n=\deg f|_H$. We show that $f$ is a generalized
polynomial of degree $n$ on $G$. Let $a_1 \stb a_{n+1}
\in G$ be arbitrary; we prove $\De _{a_1} \ldots \De _{a_{n+1}} f=0$.
We have 
\begin{equation}\label{e9}
\De _{a_1} \ldots \De _{a_{n+1}} f(x)=\sum_\va  (-1)^{|\va |+1}
f(x+\va _1 a_1 +\ldots + \va _{n+1} a_{n+1} ),
\end{equation}
where $\va =(\va _1 \stb \va _{n+1} )$ runs through all $0-1$ sequences
of length $n+1$, and $|\va |= \sum_{i=1}^{n+1} \va _i$. 
Let $x\in G$ and $\ep >0$ be fixed. Since $f$ is continuous, there is a
neighbourhood $U$ of zero such that
$$\V f(x+\va _1 a_1 +\ldots + \va _{n+1} a_{n+1} ) -f(x'+\va _1 h_1 +\ldots +
\va _{n+1} h_{n+1} )\V _E <\ep$$
for every $\va$, whenever $x'\in U+x$ and $h_i \in U+a_i$ $(i=1\stb n+1)$. 
Choosing elements $x'\in (U+x)\cap H$ and
$x_i \in (U+a_i )\cap H$ $(i=1\stb n+1)$, and noting that
$\De _{x_1} \ldots \De _{x_{n+1}} f(x')=0$ by $\deg f|_H =n$, we can see that
$\V \De _{a_1} \ldots \De _{a_{n+1}} f(x) \V _E < 2^{n+1}\ep$. Since this is true
for every $x\in G$ and $\ep >0$, it follows that $\De _{a_1} \ldots \De _{a_{n+1}}
f =0$. Thus $f$ is a generalized polynomial.

Since $f|_H$ is a polynomial, $L_{f|H}$ is of finite dimension. 
Let $T_{k_1} f|_H \stb$$T_{k_N} f|_H$ be a basis of $L_{f|H}$, and let
$V$ denote the linear hull of the functions $T_{k_1} f \stb T_{k_N} f$.
If $h\in H$, then $T_h f|_H \in L_{f|H}$, and thus there are complex numbers
$c_1 \stb c_N$ such that $T_h f|_H =\sumiN c_i T_{k_i} f|_H$; that is,
\begin{equation}\label{e10}
f(x+h)=\sumiN c_i \cd f(x+k_i )
\end{equation}
for every $x\in H$. Since $f$ is continuous and $H$ is dense in $G$, it follows
that \eqref{e10} holds for every $x\in G$. Thus $T_h f\in V$ for every $h\in H$.

Recall that $C(G,E)$, the set of continuous functions mapping $G$ into $E$,
endowed with the topology of uniform convergence on compact
sets is a topological vector space. Then $V$ is a closed subspace of $C(G,E)$,
as this is true for every finite dimensional subspace (see
\cite[Theorem 1.21]{Ru}).

We show that $L_f \su V$. Let $g\in G$ be arbitrary; we prove $T_g f\in V$.
Since $V$ is closed, it is enough to show that for every neighbourhood $W$
of $T_g f$ there is a $\phi \in V \cap W$. Let $K \su G$ be compact and $\ep >0$
be such that $\phi \in W$ whenever $\V \phi (x)-T_g f(x)\V _E <\ep$
for every $x\in K$. Since $f$ is continuous, there is a neighbourhood $Z$ of
$g$ such that $\V f (x +h )-f(x +g )\V _E <\ep$
for every $h\in Z$ and $x\in K$. Since $H$ is dense in $G$,
we can choose such a $h\in H$. Then we have $T_h f\in V\cap W$, proving
$L_f \su V$. Since $V$ is of finite dimension, so is $L_f$. \hfill $\square$
\begin{cor}\label{c1}
If $G=\rp$ with the Euclidean topology and $E$ is a Banach space, then
properties {\rm (i)-(iv)} listed in Theorem \ref{t6} are equivalent.
\end{cor}
\proof $\qq ^p$ is a dense subgroup of $\rp$ with $r_0 (\qq ^p )=p$, and thus
Theorem \ref{t7} applies. \hfill $\square$

\begin{remark} \label{r2}
{\rm We show that, in general, none of the implications \break
(i)$\akkor$(ii)$\akkor$(iii)$\akkor$(iv) can be reversed.

We saw already that the identity function defined on a Banach space of infinite
dimension is a w-polynomial but not a polynomial.

We give another example. 
Let $F$ be the free Abelian group of countable rank. We represent $F$ as
the set of sequences $x=(x_1 ,x_2 ,\ldots )$ such that $x_i$ is an integer for
every $i$ and $x_i =0$ if large enough. Let $F$ be endowed with the discrete
topology.

Let $E$ be a Banach space of infinite dimension, and let the elements
$e_1 ,e_2 ,\ldots \in E$ be linearly independent over $\cc$. Then the function
$f(x)=\sum_{i=1}^\infty x_i e_i$ $(x\in F)$ is a w-polynomial, since
$(u\circ f)(x)=\sum_{i=1}^\infty u(e_i )x_i$ is additive, hence a polynomial
on $F$ for every $u\in E^*$. On the other hand, $f$ is not a polynomial,
as $R(f)$ is of infinite dimension.

Clearly, a complex valued function is a w-polynomial if and only if it is a
polynomial. Now it is known that the function $f(x)=\sum_{i=1}^\infty x_i^2$
$(x\in F)$ is a complex valued generalized polynomial, but not a polynomial,
hence not a w-polynomial. In fact, it is easy to see
that $f$ is a  generalized polynomial of degree $2$. On the other hand, $f$
is not a polynomial, as $\dim L_f =\infty$ (see \cite{Sz3}).

Finally, 
$P(x)=\sum_{1}^\infty x_i^i$ $(x\in F)$ is a complex valued local polynomial, but
not a generalized polynomial (see \cite[Proposition 1]{L}).}
\end{remark}

\bsk
We say that a continuous function $f\in E ^G$ is a {\it w-exponential
polynomial}, if $u\circ f$ is a (complex valued) exponential polynomial for
every $u\in E^*$.
\begin{lemma}\label{l1}
If $f\in E^G$ is a w-exponential polynomial (in particular, if $f$ is a
w-polynomial), then there exists a positive integer $N$ such that
$\dim L_{u\circ f} \le N$ for every $u\in E^*$.
\end{lemma}
\proof
It is easy to see that
\begin{equation}\label{e4}
L_{u\circ f}= \{ u\circ \phi \colon \phi \in L_f\} 
\end{equation}
for every $f\in E^G$ and $u\in E^*$.
If $f\in E^G$ is a w-exponential polynomial and $u\in E^*$, then $u\circ f$ is a
complex valued exponential polynomial, and thus $\dim L_{u\circ f}$ is finite.
For every positive integer $n$, let $E_n^*$ be the set of linear functionals
$u\in E^*$ such that $\dim L_{u\circ f} <n$. Then we have $E^* =\bigcup_{n=1}^\infty
E_n^*$.

We prove that $E_n^*$ is closed. Clearly, $u\in E_n^*$ if and only if,
for every $f_1 \stb f_n \in L_f$, the complex valued functions
$u\circ f_1 \stb u\circ f_n$ are linearly dependent over $\cc$. This is true
if and only if the determinant $\det |u(f_i (x_j ))|_{i,j=1\stb n}$ is zero
for every $x_1 \stb x_n$ (see \cite[Lemma 1, p. 229]{AD}). It is easy to see
that for every
$f_1 \stb f_n \in L_f$ and $x_1 \stb x_n \in G$ the set of linear functionals
$u\in E^*$ such that $\det |u(f_i (x_j ))|_{i,j=1\stb n} =0$ is closed.
Thus $E_n^*$, the intersection of these closed sets, is also closed.

The Baire category theorem implies that there is an $n$ with
$\inter E_n^* \ne \emp$. Suppose $B(u_0 ,r)\su E_n^*$, where $B(u_0 ,r)$ is the
ball with center $u_0$ and radius $r$.

Let $u\in E^*$ be arbitrary. Then there is a $\la \in \cc \se \{ 0\}$ such that
$u_0 +\la u\in B(u_0 ,r)\su E_n^*$. By \eqref{e4}, the linear space 
$\{ u_0 \circ \phi + \la \cd u\circ \phi \colon \phi \in L_f \}$
equals $L_{(u_0 +\la u)\circ f}$, hence is of dimension $<n$.
Since the linear space $\{ u_0 \circ \phi \colon \phi \in L_f \} =
L_{u_0 \circ f}$ is also of dimension $<n$, it follows that the dimension of
$L_{u\circ f} =\{ u\circ \phi \colon \phi \in L_f \}$ is less
than $2n$. \hfill $\square$

\begin{remark}
{\rm It is easy to prove that if $f\colon G\to \cc$ is a complex valued
polynomial, then $\deg f <\dim L_f$ (see \cite[Proposition 4]{139}).
If $f\in E^G$ is a w-polynomial, then $\dim L_f$ can be infinite (see Remark
\ref{r2}).
However, by the previous lemma, there is a smallest integer $N(f)$ such that
$\dim L_{u\circ f} \le N(f)$ for every $u\in E^*$.

By Theorem \ref{t6}, every w-polynomial is a generalized polynomial, and thus
has a degree. We show that $\deg f <N(f)$ holds for every w-polynomial $f$. 

If $u\in E^*$, then $u\circ f$ is a polynomial, therefore, a generalized
polynomial. Then, by (iii) of Theorem \ref{t1}, there is an $u_0 \in E^*$
such that $\deg (u_0 \circ f)= \deg f$. 
Since $u_0 \circ f$ is a complex valued polynomial, we obtain
$$\deg f=\deg (u_0 \circ f) < \dim L_{u_0 \circ f}  \le N(f). $$

We also note that $N(f)$ is not bounded from above by any function of $\deg f$.
If, for example, $G=\rn$ and $f=x_1^2 +\ldots +x_n^2$, then $\deg f =2$. On the
other hand, $L_f$ is generated by the functions $f, x_i$ $(i=1\stb n)$ and the
constants, and thus $\dim L_f =n+2$. 
}
\end{remark}

\bsk
Our next aim is to prove the following description of w-exponential polynomials.
\begin{theorem} \label{t8}
A function $f\in E^G$ is a w-exponential polynomial if and only if there are
finitely many w-polynomials $p_1 \stb p_n \in E^G$ and complex valued
exponentials $m_1 \stb m_n \in \cc ^G$ such that $f=m_1 p_1 +\ldots +m_n p_n$.
\end{theorem}

An operator $D\colon E^G \to E^G$ is called
a {\it difference operator}, if $D$ is the linear combination with complex
coefficients of finitely many translation operators. Note that $\De _g =T_g -
T_0$ is also a difference operator. 
\begin{lemma}\label{l2}
For every finite set of distinct exponentials $\{ m_1 \stb m_n \}$ and for
every integer $s\ge -n$ there exists a finite set $\dd$ of difference operators
with the following property: whenever $p_1 \stb p_n$ are complex valued
generalized polynomials on $G$ such that $\sumin \deg p_i \le s$ and
$f=\sumin p_i \cd m_i$, then for every $1\le i\le n$ there is a $D\in \dd$ with
$p_i \cd m_i = D f$.
\end{lemma}
\proof
If $n=1$, then $\dd =\{ T_0 \}$ works, independently of $s$. Therefore, we may
assume $n>1$.

We prove by induction on $s$. Note that, by definition, the degree of the
identically zero function is $-1$. If $s=-n$, then $p_i =0$ for every
$i$, and thus $\dd =\{ {\bf 0} \}$ works, where {\bf 0} denotes the identically
zero operator (which maps every function to the identically zero function).

Suppose that $s>-n$, and that the statement is true for the smaller values.
Let $\dd$ be a finite set of difference operators such that whenever
$q_1 \stb q_n$ are generalized
polynomials with $\sumin \deg q_i \le s-1$ and $f_1 =\sumin q_i m_i$, then for
every $i$ there is a $D\in \dd$ such that $D f_1 =q_i m_i$. We may assume that
${\bf 0} \in \dd$.

From $s>-n$ it follows that $\max_{1\le i\le n} \deg p_i \ge 0$. We may assume that
$\deg p_1 \ge 0$. Since $m_1 \ne m_n$ by assumption, we can fix an element
$g\in G$ such that $m_1 (g)\ne m_n (g)$. We put $D_0 =T_g -m_1 (g)\cd T_0$. Then 
$$D_0 (p\cd m)=T_g p\cd T_g m -m_1 (g)\cd p\cd m=(m(g)\cd T_g p - m_1 (g)\cd p)
\cd m$$
for every $p\in \cc ^G$ and exponential $m$. Therefore, we have
\begin{equation}\label{e3}
  D_0 f= m_1 (g)\cd \De _g p_1 \cd m_1 +\sum_{i=2}^{n} ((m_i (g)\cd T_g p_i -
  m_1 (g)\cd p_i )\cd m_i .
\end{equation}
Since $\deg \De _g p_1 <\deg p_1$ and $\deg ( (m_i (g)\cd T_g p_i -  m_1 (g)\cd
p_i )\le \deg p_i$, it follows from the choice of $\dd$ that
$m_1 (g)\cd \De _g p_1 \cd m_1 =D_{1} D_0 f$ for some $D_1 \in \dd$. Similarly,
for every $i$ such that $p_i \ne 0$ there is a $D_i \in \dd$ with
$m_i (g) \De _g p_i \cd m_i =D_{i} D_0 f$. If $p_i =0$ then we can take
$D_i ={\bf 0}$, and thus there is such a $D_i$ in both cases. By \eqref{e3}
and by the choice of $\dd$ we have
$$(m_n (g)\cd T_g p_n -  m_1 (g)\cd p_n )\cd m_n =E_n D_0 f$$
with a suitable $E_n \in \dd$. Since
\begin{align*}
(m_n (g)\cd T_g p_n  &-  m_1 (g)\cd p_n )\cd m_n =\\
& =m_n (g)\cd \De _g p_n \cd m_n   + (m_n (g) -m_1 (g))\cd p_n \cd m_n 
\end{align*}  
and $m_n (g)\cd \De _g p_n \cd m_n =D_{n} D_0 f$, we obtain
$$E_n D_0 f= D_n D_0 f +(m_n (g) -m_1 (g))\cd p_n \cd m_n$$
and 
$$p_n \cd m_n =c\cd E_n D_0 f -c\cd D_n D_0 f ,$$
where $c=1/(m_n (g) -m_1 (g))$. Therefore, if we add the operators $c D D_0 -
cD' D_0$ $(D, D'\in \dd )$ to $\dd$, then $p_n \cd m_n =E f$ will hold for
a suitable $E$ belonging to the enlarged $\dd$.
(Note that the element $g$ does not depend on the functions $p_1 \stb p_n$, only
on $m_1$ and $m_n$.) The same argument provides finitely many operators
such that if we add them to $\dd$ then, for every $i=1\stb n$, $p_i m_i =Ef$
will hold for a suitable $E$ belonging to the enlarged $\dd$. \hfill $\square$

\bsk
{\bf Proof of Theorem \ref{t8}.} 
Suppose $f=m_1 p_1 +\ldots +m_n p_n$, where $p_1 \stb p_n$ are w-polynomials and
$m_1 \stb m_n$ are complex valued exponentials. If $u\in E^*$, then
$u\circ f=m_1 \cd u\circ p_1 +\ldots +m_n \cd u\circ p_n$. Since $p_i$ is a
w-polynomial, it follows that $u\circ p_i$ is a complex valued polynomial for
every $i=1\stb n$, and thus $u\circ f$ is an exponential polynomial.
This is true for every $u\in E^*$, proving that $f$ is a w-exponential
polynomial.

Now suppose that $f$ is a w-exponential polynomial.
By Lemma \ref{l1}, there is a positive integer $K$ such that $\dim L_{u\circ f}
<K$ for every $u\in E^*$.

Let $\pp$ denote the set of all functions $p\cd m$ such that 
$p\colon G\to \cc$ is a polynomial and $m\colon G\to \cc$ is an exponential.
If $u\in E^*$, then $u\circ f$ is an exponential polynomial, and thus
it is the sum of finitely many elements of $\pp$. In other words,
for every $u\in E^*$ there exists a finite set $\pp _u \su \pp$ such that
$u\circ f=\sum_{p\cd m\in \pp _u} p\cd m$. 

Let $\mmm$ denote the set of those exponentials $m$ for which there exist $u\in
E^*$ and a nonzero polynomial $p$ such that $p\cd m\in \pp _u$.
We prove that $\mmm$ contains less than $K$ distinct exponentials.

Suppose this is not true, and let $m_1 \stb m_{K}$ be distinct exponentials
in $\mmm$. We may assume that for every $u\in E^*$ and $1\le i\le K$ there is a
unique polynomial $p_{u,i}$ such that $p_{u,i} \cd m_i \in \pp _u$. Indeed, if
$\pp _u$ does not contain such a product, then we add $0\cd m_i$ to $\pp _u$,
and put $p_{u,i} =0$.

For every $1\le i\le K$ we have $m_i \in \mmm$, and thus there is an $u_i \in
E^*$ such that $p_{u_i ,i} \ne 0$. We show that there are complex numbers
$\la _1 \stb \la _{K}$ such that
\begin{equation}\label{e5}
\sum_{i=1}^K \la _i  p_{u_i ,j} \ne 0
\end{equation}
for every $j=1\stb K$. Indeed, for a fixed $j$, the set of $K$-tuples
$(\la _1 \stb \la _{K} )$ such that $\sum_{i=1}^K \la _i  p_{u_i ,j} =0$
is a linear subspace $L_j$ of $\cc ^K$. Since $p_{u_j ,j} \ne 0$,
the subspace $L_j$ does not contain the vector $(0\stb 0,1,0\stb 0)$ having
$1$ as the $j$th coordinate. Therefore, $L_j$ is a proper subspace of $\cc ^K$.
Now $\cc ^K$ is not the union of finitely many proper subspaces, therefore,
we must have \eqref{e5} for every $j=1\stb K$ with a suitable 
$(\la _1 \stb \la _{K} )$.

Let $u=\sum_{i=1}^K \la _i u_i$. Then $p_{u ,j} =\sum_{i=1}^K \la _i  p_{u_i ,j}
\ne 0$ for every $j=1\stb K$. That is, in the representation of $u\circ f$
as a sum of functions $p\cd m\in \pp$, each of $m_1 \stb m_K$ appears with
a nonzero polynomial factor.

Now we need the following result: if $V$ is a translation invariant 
linear subspace of $\cc ^G$ and $\sumin p_i m_i \in V$, where $p_i \in \cc^G$
is a nonzero generalized polynomial for every $i=1\stb n$ and $m_1 \stb m_n$ are
distinct complex valued exponentials, then $m_i \in V$ for every $i=1\stb n$.
This is proved, e.g., in \cite[Lemma 6]{L} in the case when $G$ is an Abelian
group. One can easily check that the same proof works in Abelian semigroups.

Since $L_{u\circ f}$ is a translation invariant 
linear subspace of $\cc ^G$ and $u\circ f=\sum _{p m \in \pp _u} p\cd m$, it
follows that $m_1 \stb m_K \in L_{u\circ f}$. Now $m_1 \stb m_K$ are linearly
independent over $\cc$, since, if $\sum_{i=1}^K c_i m_i =0$, where
$c_1 \stb c_K \in \cc$, then the unique representation of the zero function
implies $c_1 =c_2 =\ldots =c_K =0$. We find that $\dim L_{u\circ f} \ge K$,
which contradicts the choice of $K$.

This contradiction proves that $\mmm$ contains less than $K$ exponentials.
Let $\mmm =\{ m_1 \stb m_n \}$, where $n<K$. Then for every $u\in E^*$
there are polynomials $p_{u ,i}$ $(i=1\stb n)$ such that
$u\circ f=\sumin p_{u ,i} \cd m_i$.

Let $s(u)=\max_{1\le i\le n} \deg p_{u ,i}$, and put $U_s =\{ u\in E^* \colon
s(u)\le s\}$
for every positive integer $s$. By the Baire category theorem we can find
an integer $s$ such that $U_s$ is of second category in $E^*$.

By Lemma \ref{l2}, there is a finite set $\dd$ of difference operators
such that, for every $u\in U_s$ and $i=1\stb n$, 
$p_{u ,i} \cd m_i =D_i (u\circ f)$ for some $D_i \in \dd$. Since $\dd$ is finite,
there is an $n$-tuple $(D_1 \stb D_n )$ such that the set $A$ of linear 
functionals $u\in U_s$ such that
\begin{equation}\label{e6}
p_{u ,i} \cd m_i =D_i (u\circ f) \qquad (i=1\stb n)
\end{equation}
is also of second category in $E^*$.
We fix such an $n$-tuple $(D_1 \stb D_n )$, and define $B$ as the set of those
functionals $u\in E^*$ for which \eqref{e6} holds. Then $A\su B$, and thus $B$
is of second category in $E^*$.

We show that $B$ is a closed linear subspace of $E^*$. Suppose $u_1 ,u_2\in B$,
$\la _1  ,\la _2 \in \cc$, and put $u=\la _1 u_1 +\la _2 u_2$. Since 
$u_j \circ f =\sumin p_{u_j ,i} \cd m_i$ $(j=1,2)$, we have
$$u\circ f= \sumin (\la _1 p_{u_1 ,i} +\la _2 p_{u_2 ,i} )\cd m_i .$$
The uniqueness of the representation gives
$$p_{u ,i} =\la _1 p_{u_1 ,i} +\la _2 p_{u_2 ,i}$$
and
\begin{align*}
D_i (u\circ f) &=\la _1 D_i (u_1 \circ f) +\la _2 D_i (u_2 \circ f) =\\
&=\la _1 p_{u_1 ,i} \cd m_i  +\la _2 p_{u_2 ,i} \cd m_i  =\\
&=p_{u ,i} \cd m_i .
\end{align*}
Thus $u\in B$, proving that $B$ is a linear subspace of $E^*$.
Let $u\in E^*$ be in the closure of $B$. Then there is a sequence of linear
functionals $u_\nu \in B$ such that $\V u_\nu - u \V \to 0$. Then
$$p_{u_\nu  ,i} \cd m_i =D_i (u_\nu \circ f)\to D_i (u\circ f)$$
as $\nu \to \infty$, for every $i=1\stb n$. Thus $p_{u_\nu  ,i} \to q_i$ pointwise,
where $q_i =D_i (u\circ f)/m_i$.

Now $u_\nu \in B\su A\su U_s$, and thus $p_{u_\nu  ,i}$ is a generalized
polynomial of degree $\le s$. It is easy to check that this property is
preserved under pointwise convergence in the discrete topology of $G$.
Therefore, $q_i$ is generalized polynomial of degree $\le s$ for every
$i=1\stb n$ in the discrete topology of $G$. Now we have
$$\sumin p_{u  ,i} \cd m_i = u\circ f= \lim_{\nu \to \infty} u_\nu \circ f = \lim_{\nu \to \infty}
\sumin p_{u_\nu  ,i} \cd m_i =\sumin q_i \cd m_i .$$
Then the uniqueness of the representation gives $q_i =p_{u  ,i}$ for every $i=1
\stb n$. Thus $p_{u  ,i} =D_i (u\circ f)/m_i$, $D_i (u\circ f) =p_{u  ,i}\cd m_i$
for every $i=1\stb n$; that is $u\in B$.

Thus $B$ is a closed subspace of $E^*$. Since $B$ is of second category, we
have $B=E^*$. Therefore, \eqref{e6} holds for every $u\in E^*$.

Let $p_i =(D_i f )/m_i$ $(i=1\stb n)$ and $\ol f =\sumin p_i \cd m_i$.
Since
$$u\circ p_i =u\circ (D_i f )/m_i =D_i (u\circ f)/m_i =p_{u  ,i}$$
is a polynomial for every $u\in E^*$, it follows that $p_i$ is a w-polynomial
for every $i=1\stb n$. Now
\begin{align*}
u\circ \ol f & =\sumin (u\circ p_i ) \cd m_i =\sumin u\circ D_i f=\\
&=\sumin D_i (u\circ f)=\sumin p_{u  ,i} \cd m_i =\\
&=u\circ f
\end{align*}
for every $u\in E^*$. Thus $\ol f =f$, which completes
the proof. \hfill $\square$.

The following result is an immediate consequence of Theorems \ref{t5}, \ref{t7}
and \ref{t8}.
\begin{cor}\label{c2}
Let $G$ be a topological Abelian group, and suppose that there is a dense
subgroup $H$ of $G$ such that $r_0 (H) <\infty$. Then a function
$f\in E^G$ is a w-exponential polynomial if and only if $f$ is an
exponential polynomial.

In particular, if $G=\rp$ with the Euclidean topology and $E$ is a Banach space,
then a function $f\in E^G$ is a w-exponential polynomial if and only if $f$ is
an exponential polynomial.
\end{cor}

\section{Vector valued harmonic analysis and synthesis on discrete Abelian
groups}

Let $G$ be a topological Abelian group, and $E$ be a Banach space.
We denote by $C(G,E)$ the set of continuous functions $f\colon G\to E$.
We equip $C(G,E)$ with the topology of uniform convergence on compact sets.
In this topology a set $U\su C(G,E)$ is open if, for every $f\in U$, there
exists a compact set $K \su G$ and there is an $\ep >0$ such that
$\phi \in U$ whenever $\phi \in C(G,E)$ is such that $\V \phi (x) -f(x) \V _E
<\ep$ for every $x\in K$. This topology makes $C(G,E)$ a locally convex
topological vector space over the complex field.

Translation invariant closed linear subspaces of $C(G,E)$ are called {\it
varieties}. If $f\in C(G,E)$, then $V_f$ denotes the smallest variety containing
$f$. Clearly, $V_f$ equals the closure of $L_f$.
We say that {\it spectral synthesis holds in} $C(G,E)$ if every variety $V$
in $C(G,E)$ is the closed linear hull of the set of exponential polynomials
contained in $V$.

It is known that spectral synthesis holds in $C(G,E)$ for every compact
Abelian group $G$ and for every Banach space $E$. In fact, this result follows
from the approximation theorem of almost periodic functions mapping
a group into a Banach space \cite{BN}. More precisely, if $G$ is a compact
Abelian group, then (i) every continuous function $f\colon G\to E$ is almost
periodic, (ii) the terms of the Fourier series of $f$ are constant multiples of
characters belonging to $V_f$, and (iii) it follows from the approximation
theorem that linear combinations of these characters approximate $f$ uniformly.
We note, however, that this special case can be proved without the machinery
of almost periodic functions, and in the next section we provide a simple
self-contained proof.

As it turns out, the situation for locally compact Abelian groups is different.
In the rest of this section we concentrate on discrete Abelian groups.
In these groups we have $C(G,E)=E^G$.
As we will see shortly, spectral synthesis does not hold in $E^G$ if $G$ is
infinite and $E$ is of infinite dimension. Actually, the situation is even
worse.

By a generalized (resp. local) exponential polynomial we mean a function of the
form $\sumin m_i \cd p_i$, where $m_i$ is an exponential and $p_i \in E^G$ is a
generalized (resp. local) polynomial on $G$ for every $i=1\stb n$. 
We say that {\it generalized (resp. local) spectral synthesis holds in} $E^G$
if every variety $V$ in $E^G$ is the closed linear hull of the set of
generalized (resp. local) exponential polynomials contained in $V$. Clearly, the
condition that  spectral synthesis holds in $E^G$ implies that
generalized spectral synthesis holds in $E^G$, and this condition, in turn,
implies that local spectral synthesis holds in $E^G$.

Let $f=\sumin m_i \cd p_i$, where $m_1 \stb m_n$ are distinct exponentials and
$p_1 \stb p_n$ are nonzero local polynomials on $G$.
One can prove that $m_i \cd p_i \in V_f$ for every $i$, moreover, there are
nonzero elements $e_1 \stb e_n \in E$ such that
$m_i \cd e_i \in V_f$ for every $i$. (See \cite[Lemma 7]{L}, where the complex
valued case is proved. One can easily check that the proof works in the general
case as well.) This result shows that if local spectral synthesis holds in
$E^G$, then {\it spectral analysis holds in} $E^G$; that is, every nonzero
variety in $E^G$ contains a function of the form $m \cd e$, where $m$ is an
exponential and $e\in E$, $e\ne 0$.
\begin{theorem} \label{t10}
If $G$ is an infinite discrete Abelian group and $E$ is a Banach space of
infinite dimension, then spectral analysis does not hold in $E^G$.
\end{theorem}
\proof
It is easy to see that if spectral analysis holds in $E^G$, then
the same is true in $C(H,E)$ for every subgroup $H$ of $G$. (For the complex
case see \cite[Lemma 4]{LSz1}. The proof in the general case is the same.)
Therefore, in order to prove the theorem, it is enough to find a subgroup
$H$ of $G$ such that spectral analysis does not hold in $C(H,E)$.

If $G$ contains an element $h$ of infinite order, then we let $H$ be the
cyclic group generated by $h$. If $G$ is torsion, then we choose a countably
infinite subset $A\su G$, and let $H$ be the subgroup generated by $A$.
Then $H$ is countably infinite, and is either cyclic, or torsion.

Let $g_1 ,g_2 ,\ldots$ be an enumeration of the elements of $H$.
If $H$ is cyclic generated by the element $h$, then we choose an enumeration
such that $g_{2n} =h^n$ $(n=1,2,\ldots )$. If $H$ is torsion, then the
enumeration can be arbitrary.

Since $E$ is of infinite dimension, it contains a basic sequence $(x_n )$
(see \cite[Corollary 3, p. 39]{D}). We may assume that $\V x_n \V =n!$
for every $n=1,2,\ldots$.

We define $f(g_n )=x_n$ for every $n=1,2,\ldots$, and prove that $V_f$
does not contain any function of the form $m \cd e$, where $m$ is an
exponential and $e\in E$, $e\ne 0$.

Suppose this is false, and let $m \cd e \in V_f$, where $m$ and $e$ are as
above. Since $H$ is countable and $m\cd e \in V_f =\cl L_f$, it follows that 
there is a sequence of functions $f_k \in L_f$ such that $f_k \to m\cd e$
pointwise on $H$. Now each $f_k$ is a linear combination of the functions
$T_{g_n} f$, and thus there is a sequence of positive integers $s_1 <s_2 <\ldots$
such that $f_k$ is a linear combination of the functions
$T_{g_n} f$ $(n=1\stb s_k )$. Let
$$f_k =\sum_{n=1}^{s_k} c_{kn} T_{g_n} f \qquad (k=1,2,\ldots ).$$
Let $e=\sum_{n=1}^\infty \al _n x_n$ with suitable complex coefficients $\al _n$.
Since the series converges in norm, it follows that $\V \al _n x_n \V \to 0$,
that is, $n!|\al _n |\to 0$. We have
$$\sum_{n=1}^{s_k} c_{kn} x_n = \sum_{n=1}^{s_k} c_{kn} f(g_n ) =f_k (0)\to m(0)\cd e=
e= \sum_{n=1}^\infty \al _n x_n$$
as $k\to \infty$. Since the coefficient functionals are continuous (see
\cite[p. 32]{D}), it follows that $\lim_{k\to \infty} c_{kn}=\al _n$ for
every $n$. Let the indices $i$ and $j$ given, and put $g=g_j -g_i$, Then
$$f_k (g)=\sum_{n=1}^{s_k} c_{kn} f(g+g_n )\to m(g)\cd e =\sum_{n=1}^\infty
(m(g)\al _n )x_n$$
as $k\to \infty$. The elements $g+g_n$ $(n=1\stb s_k )$ are distinct, and then
so are $f(g+g_n )$. If $k$ is large enough, then $i<s_k$, $f(g+g_i )=f(g_j )=
x_j$, and thus $f_k (g)$ is a linear combination of finitely many of the
elements $x_n$ including $x_j$. Since the coefficient functionals are
continuous, it follows that the coefficient of $x_j$ converges to $m(g)\al _j$
as $k\to \infty$. That is, we have $\lim_{k\to \infty} c_{ki} = m(g)\al _j$.
However, as $c_{ki} \to \al _i$, we find
$$\al _i =m(g)\al _j =m(g_j -g_i )\al _j =(m(g_j )/m(g_i ))\al _j$$
and $\al _i m(g_i )=\al _j m(g_j )$. This is true for every $i$ and $j$,
and thus there is a complex number $c$ such that $\al _i m(g_i )=c$ for every
$i=1,2,\ldots$.

Since $e\ne 0$, we have $\al _n \ne 0$ for at least one $n$, and thus $c\ne 0$.
Now we consider the two cases concerning the group structure of $H$.
If $H$ is torsion, then the value of $m(g)$ is a root of unity for every $g\in H$. Then $|\al _n |= |c/m(g_n )|=|c|$ for every $n$. 
This, however, is impossible by $n!|\al _n |\to 0$.

Next suppose that $H$ is cyclic with generator $h$. Then we have 
$g_{2n}=h^n$ for every $n=1,2,\ldots$. Let $m(h)=\la$, then
$m(g_{2n} )=\la ^{n}$ and $|\al _{2n} |= |c/m(g_{2n} )|=|c|\cd |\la |^{-n}$ for
every $n$. This, again, contradicts $n!|\al _n |\to 0$, completing the
proof. \hfill $\square$

\bsk
Returning to the question of spectral synthesis, the previous result shows
that spectral synthesis can hold in $E^G$ only if $G$ is finite or
$E$ is of finite dimension. If $G$ is finite, then every element $f\in E^G$
is an exponential polynomial, as $L_f$ is of finite dimension. Therefore, in
this case spectral synthesis does hold.

If $E$ is of finite dimension, then spectral synthesis in $E^G$ is still 
not automatic. Indeed, if spectral synthesis holds in $E^G$, then it
also holds in $C(G,\cc) =\cc ^G$. Now it is known that spectral synthesis holds
in $C(G,\cc )$ if and only if $r_0 (G)$ (the torsion free rank of $G$) is finite
(see \cite[Theorem 1]{LSz2}). So the only cases left are when
$r_0 (G)$ is finite and $E$ is of finite dimension. In the next theorem
we show that spectral synthesis does hold in these cases.

We also consider local spectral synthesis and spectral analysis.
Note that spectral analysis holds in $C(G,\cc )$ 
if and only if $r_0 (G)$ is less than continuum (see \cite[Theorem 1]{LSz1}).
Also, there exists an uncountable cardinal $\kappa$ such that
local spectral synthesis holds in $C(G,\cc )$ if and only if $r_0 (G)<\kappa$
(see \cite[Theorem 3]{L}).
In particular, local spectral synthesis holds in $C(G,\cc )$ for every
countable discrete Abelian group.
\begin{theorem} \label{t11}
Let $G$ be a discrete Abelian group, and let $k$ be a positive integer.
\begin{enumerate}[{\rm (i)}]
\item  If $r_0 (G)$ is finite, then spectral synthesis holds in
$C(G, \cc ^k )$.
\item  If $r_0 (G) <\kappa$, then 
local spectral synthesis holds in $C(G, \cc ^k )$.
\end{enumerate}
\end{theorem}

For every set $V$ of maps
$f\colon G\to \cc ^k$ we shall denote by $\ol V$ the set of maps
\begin{equation}\label{e1}
\zk \times G \ni (t_1 \stb t_k ,x) \mapsto t_1 f_1 (x)+\ldots +t_k f_k (x)+g(x),
\end{equation}
where $(f_1 \stb f_k )\in V$ and $g\colon G\to \cc$. 

\begin{lemma}\label{l3}
If $V$ is a variety of maps $f\colon G\to \cc ^k$,
then $\ol V$ is a variety on $\zk \times G$.
\end{lemma}

\proof
Suppose that $V$ is a variety. It is clear that $\ol V$ is a translation invariant linear space. We show that $\ol V$ is closed.

Let $e_i =(\de _{1i}\stb \de_{ki},0) \in \zk \times G$ for every $i=1\stb k$,
where $\de _{ji}$ is the Kronecker delta. It is clear that if the function $F$
is defined by \eqref{e1}, then $\De _{e_i} F(t_1 \stb t_k ,x)=f_i (x)$ for every
$i=1\stb k$ and $(t_1 \stb t_k ,x) \in \zk \times G$. 
Suppose that $h\colon (\zk \times G) \to \cc$ is in the closure
of $\ol V$. From the previous observation it follows that
$\De _{e_i} h$ does not depend on the variables $t_1 \stb t_k$. That is,
there are functions
$h_1 \stb h_k \colon G\to \cc$ such that
$\De _{e_i} h(t_1 \stb t_k ,x)=h_i (x)$ for every $i=1\stb k$ and $(t_1 \stb t_k ,x) \in \zk \times G$.
Let
$$s(t_1 \stb t_k ,x)=t_1 h_1 (x)+\ldots +t_k h_k (x)$$
for every $(t_1 \stb t_k ,x) \in \zk \times G$, and put $g=h-s$. Then
$\De _{e_i} g=0$ for every $i=1\stb k$. Thus $g=h-s$  does not depend on the variables $t_1 \stb t_k$. Therefore, we have
$$h(t_1 \stb t_k ,x)=t_1 h_1 (x)+\ldots +t_k h_k (x) +g(x)$$
for every $(t_1 \stb t_k ,x) \in \zk \times G$.
We prove that $(h_1 \stb h_k )\in V$. Since $V$ is a variety, it is enough to show that $(h_1 \stb h_k )$ is in the closure of $V$.

Let the finite set $X\su G$ and the positive number $\ep$ be given.
Since $h$ is in the closure of $\ol V$, there is a function $f\in \ol V$ such that
it is closer to $h$ than $\ep /2$ at each point $(t_1 \stb t_k ,x)$,
where $t_i =0,1$ for every $i=1\stb k$ and $x\in X$. Let
$$f(t_1 \stb t_k ,x)=t_1 f_1 (x)+\ldots +t_k f_k (x) +g_1 (x),$$
where $(f_1 \stb f_k )\in V$.
Since $\De _{e_i} h(0\stb 0,x)=h_i (x)$ and $\De _{e_i} f(0\stb 0 ,x)=f_i (x)$ for
every $x\in X$, it follows that $|h_i (x)-f_i (x)|<\ep$ for every $x\in X$
and $i=1\stb k$. This proves that $(h_1 \stb h_k )$ is in the closure of $V$.
Thus $(h_1 \stb h_k ) \in V$ and $h\in \ol V$, showing that $\ol V$ is a
variety. \hfill $\square$

\bsk \noi
{\bf Proof of Theorem \ref{t11}}. (i) Let $V$ be a variety of maps $f\colon
G\to \cc ^k$. By Theorem \ref{t3}, a function $f\colon G\to \cc ^k$ is an
exponential polynomial if and only if $f=(f_1 \stb f_k )$, where
$f_1 \stb f_k$ are complex valued exponential polynomials.

We have to show that if $r_0 (G)$ is finite, then the set of maps
$$\{ (p_1 \stb p_k ) \in V \colon p_1 \stb p_k \ \text{are exponential
polynomials}\}$$
is dense in $V$. Let $(f_1 \stb f_k ) \in V$, and let the finite set $X\su G$
and the positive
number $\ep$ be given. If $r_0 (G)$ is finite, then $r_0 ( \zk \times G)$
is also finite. Then, by \cite[Theorem 1]{LSz2}, spectral synthesis
holds in every variety on $\zk \times G$. By Lemma \ref{l3}, $\ol V$ is a
variety. Since the function $f=t_1 f_1 (x)+\ldots +t_k f_k (x)$ belongs to
$\ol V$, it follows that there is an exponential polynomial $p\in \ol V$
such that $p$ is closer to $f$ than $\ep /2$ at each point $(t_1 \stb t_k ,x)$,
where $t_i =0,1$ for every $i=1\stb k$ and $x\in X$. Let
$$p(t_1 \stb t_k ,x)=t_1 p_1 (x)+\ldots +t_k p_k (x) +g_2 (x).$$
Then $(p_1 \stb p_k )\in V$ by $p\in \ol V$, and $|p_i (x)-f_i (x)|<\ep$ for
every $x\in X$ and $i=1\stb k$. Since $p$ is an exponential polynomial, so is
$p_i =\De _{e_i} p(0\stb 0 ,x)$ for every $i$.
This proves that the maps $(p_1 \stb p_k ) \in V$, where $p_1 \stb p_k$ are exponential polynomials constitute a dense subset of $V$. This proves (i).

\msk \noi
The proof of (ii) is similar to that of (i). If $r_0 (G) <\kappa$, then
$r_0 ( \zk \times G) <\kappa$. Thus local spectral synthesis
holds in the variety $\ol V$. 
If $f=t_1 f_1 (x)+\ldots +t_k f_k (x) \in \ol V$, there is a
local exponential polynomial $p\in \ol V$
such that $p$ is closer to $f$ than $\ep /2$ at each point $(t_1 \stb t_k ,x)$,
where $t_i =0,1$ for every $i=1\stb k$ and $x\in X$. As we saw above, this implies that $|p_i (x)-f_i (x)|<\ep$ for
every $x\in X$ and $i=1\stb k$. Since $p$ is a local exponential polynomial,
so is $p_i =\De _{e_i} p(0\stb 0 ,x)$ for every $i$.
This proves that the maps $(p_1 \stb p_k ) \in V$, where $p_1 \stb p_k$ are
local exponential polynomials constitute a dense subset of $V$.
\hfill $\square$

\begin{cor}\label{c3}
If $r_0 (G)<\kappa$ (in particular, if $r_0 (G)$ is countable), then
spectral analysis holds in $C(G, \cc ^k )$.
\end{cor}
\proof
Let $V\su C(G, \cc ^k )$ be a nonzero variety.
By (ii) of Theorem \ref{t11}, local spectral synthesis holds 
in $C(G, \cc ^k )$, and thus there is a nonzero local exponential polynomial
$f\in V$. Let $f=\sumin m_i \cd p_i$, where $m_i$ is an exponential and $p_i \in
E^G$ is a local polynomial for every $i=1\stb n$. We may assume that
$m_1 \stb m_n$ are distinct and $p_1 \stb p_n$ are nonzero. Then,
by \cite[Lemma 7]{L}, $m_i \cd e_i \in V$ for every $i=1\stb n$ with
nonzero $e_1 \stb e_n \in E$. In fact, \cite{L}
deals with complex valued maps, but the argument of the proof of
\cite[Lemma 7]{L} works for vector valued functions as well. This proves
that spectral synthesis holds in $V$. \hfill $\square$

Note that if $k=1$, then $\kappa$ can be replaced by $2^\om$ in the statement
of Corollary \ref{c3} (see \cite{LSz1}). Since $\om _1 \le \kappa \le 2^\om$,
it makes no difference under the continuum hypothesis. Still, it would be
interesting to see if $\kappa$ can be replaced by $2^\om$ in the cases $k>1$
as well.

\section{Appendix: a proof of spectral synthesis on compact Abelian groups.}
\begin{theorem} \label{t9}
If $G$ is a compact Abelian group and $E$ is a Banach space, then spectral
synthesis holds in $C(G,E)$.
\end{theorem}
\proof
We denote by $M(G)$ the set of complex valued regular Borel measures on $G$
with finite total variation $\V \mu \V$.
If $f\in C(G,E)$ and $\mu \in M(G)$, then we define $\mu *f(x)$ as
$\int_G f(x-t)\, d\mu (t)$. The integral makes sense for every $x\in G$, as we
integrate a continuous function mapping a compact Hausdorff space into a Banach
space (see \cite[Theorem 3.27]{Ru}).

We show that $\mu *f\in V_f$ for every $f\in C(G,E)$ and $\mu \in M(G)$. 
Let $\ep >0$ be given. By the uniform continuity of $f$ we can find a
neighbourhood $V$ of $0$ such that $\V f(x)-f(y)\V _E <\ep /\V \mu \V$ whenever
$x-y\in V$. Let $U$ be a symmetric neighbourhood of $0$ such that $U+U\su V$.

Since $G$ is compact, there are points
$x_i \in G$ $(i=1\stb n)$ such that $G=\bigcup_{i=1}^n (U+x_i )$.
Choose a partition $\{ E_1 \stb E_n  \}$
of $G$ into finitely many Borel sets such that $E_i \su U+x_i$ for every
$i=1\stb n$. Deleting the empty sets from the partition, we may assume that
$E_i \ne \emp$ for every $i$. Choose a point $t_i \in E_i$ for every
$i=1\stb n$, and put
\begin{equation}\label{e11}
z(x)=\int_G f(x-t)\, d\mu (t) -  \sumin \mu (E_i )f(x-t_i ) \qquad (x\in G).
\end{equation}
Then $\V z(x)\V _E \le \ep$ for every $x\in G$. Indeed, suppose $\V z(x)\V _E
> \ep$ for some $x$. Then there is a
$u\in E^*$ such that $\V u\V \le 1$ and $|u(z(x))|> \ep$. Now
\begin{align*}
u(z(x)) &=\int_G (u\circ f)(x-t)\, d\mu (t) - \sumin \mu (E_i )\cd
(u\circ f)(x-t_i )=\\
&= \sumin \int_{E_i} \left( (u\circ f)(x-t) - (u\circ f)(x-t_i ) \right) \,
d\mu (t).
\end{align*}
However, since $t-t_i \in V$ and $\V u(f(x-t))-u(f(x-t_i ))\V _E <
\ep /\V \mu \V$ whenever $t\in E_i$, we have $|u(z(x))|\le \ep$, which is
impossible.

Now the function $\sumin \mu (E_i )f(x-t_i )$ is in $L_f$. By \eqref{e11} 
we find that $\mu *f$ can be uniformly approximated by functions from $L_f$,
and thus $\mu *f\in V_f$.

If $f\in C(G,E)$ and $g\in C(G,\cc )$ then $f*g(x)$ is defined by
$\int_G f(x-t)g(t) \, dt$, where $dx$ is the Haar measure. Clearly,
$f*g=\mu _g *f$, where $\mu _g (H)=\int_H g\, dt$ for every Borel $H\su G$.
Thus $f*g\in V_f$.

It is clear that $C(G,E)$ is a Banach space with the norm $\V f\V =\sup_{x\in G}
\V f \V _E$. If $\La \in C(G,E)^*$ and $f\in C(G,E)$, then we put
$\La *f (x)= \La (g_x )$, where $g_x$ is the function $t\mapsto f(x-t)$
$(t\in G)$.
Using the uniform continuity of $f$ it is easy to see that $\La *f \in
C(G,\cc )$. Our next aim is to show that if $\La \in C(G,E)^*$, $f\in C(G,E)$
and $g\in C(G,\cc )$, then 
\begin{equation}\label{e12}
(\La *f )*g=\La *(f*g).
\end{equation}
First we assume that the linear span $N(f)$ of $R(f)$ is of finite dimension.
Let $e_1 \stb e_n$ be a basis of $N(f)$. Then there are functions $f_i
\colon G\to \cc$ such that $f=\sumin f_i \cd e_i$. 
There are linear functionals $u_1 \stb u_k \in E^*$ such that $u_i (e_i)=1$
and $u_i (e_j )=0$ for every $1\le i,j\le k$, $i\ne j$. Thus $f_i =u_i \circ f$
for every $i=1\stb k$. Since $f$ is continuous, we can see that so are $f_1
\stb f_k$.

Let $\La _i (g)=\La (g\cd e_i )$ for every $g\in C(G,\cc )$ and $i=1\stb n$.
It is clear that $\La _i$ is a bounded linear functional of $C(G,\cc )$.
By the Riesz representation theorem, there are measures $\mu _1 \stb \mu _n$
such that $\La _i (g)=\int _G g \, d\mu _i$ for every $g\in C(G,\cc )$ and
$i=1\stb n$. Then $\La _i *g=\mu _i *g$ for every $g\in C(G,\cc )$ and
$i=1\stb n$, and thus
\begin{align*}
  (\La *f )*g&= \left( \La *\sumin f_i \cd e_i \right) *g=
   \left( \sumin  \La * (f_i \cd e_i )\right) *g=\\
&= \left( \sumin  \La _i *f_i  \right) *g=
\left( \sumin  \mu_i *f_i  \right) *g=\\
  &=\sumin (\mu_i *f_i ) *g= \sumin \mu_i *(f_i  *g)=\\
  &=\sumin \La _i *(f_i  *g)=
\sumin  \La  *(( f_i \cd e_i )  *g)=\\
&=\La *(f*g).
\end{align*}
If $f\in C(G,E)$ is arbitrary, then a standard partition-of-unity
argument\footnote{Let $\ep >0$ be given. Let $V$ be a
neighbourhood of $0$ such that $\V f(x)-f(y)\V _E <\ep$ whenever
$x-y\in V$. Let $U$ be another neighbourhood of $0$ such that $\cl U\su V$,
and let $\phi \colon G\to \rr$ be a continuous function such that
$\phi >0$ on $U$, $0\le \phi \le 1$ on $V$ and $\phi =0$ on $G\se V$. 
Let $x_i \in G$ $(i=1\stb n)$ be such that $G=\bigcup_{i=1}^n (U+x_i )$,
and put $\Phi (x) =\sumin \phi (x-x_i )$. Then $\Phi$ is a positive continuous
function on $G$. Put $f_\ep =\Phi \am \cd \sumin \phi (x-x_i )\cd f(x_i )$.
Then $f_\ep \in C(G,E)$, and the linear span of $R(f_\ep )$ is of finite
dimension. Now we have
$f(x)-f_\ep (x) = \Ph \am \sumin \phi (x-x_i )\cd (f(x)-f(x_i ))$
for every $x$. If $x\notin V+x_i$, then $\phi (x-x_i )=0$. If
$x\in V+x_i$, then $0\le \phi (x-x_i )\le 1$ and $\V f(x)-f(x_i ) \V _E <\ep$,
and thus $\V f(x)-f_\ep (x)\V _E <\ep$.}
shows that $f$ can be uniformly approximated by functions $f_n$ such that
the linear span of $R(f_n )$ is of finite dimension for every $n$.
Suppose $\V f-f_n \V _E <1/n$ $(n=1,2,\ldots )$. Then $\V f*g-f_n *g\V _E \le
\V g \V _1 /n$ on $G$,
$$|\La *(f*g)-\La * (f_n *g)|\le \V \La \V \cd \V g \V _1 /n ,$$
$$\left| \La *f -\La *f _n \right| =\left| \La * (f -f _n ) \right| 
\le \V \La \V /n,$$
$$\left\V (\La *f)*g -(\La *f _n )*g \right\V \le \V \La \V \cd \V g \V _1 /n
.$$
Since $(\La *f_n )*g=\La *(f_n *g)$, it follows that the norm of the difference
of the two sides of \eqref{e12} is at most $c/n$ for every $n$, where $c$
only depends on $\La$ and $g$. This proves \eqref{e12}.

Let $\Ga$ denote the dual of $G$. If $\ga \in \Ga$, then the invariance of the
Haar measure gives
$$(f*\ga )(x)= \int_G \ga (x-t)f(t) \, dt =e_\ga \cd \ga (x),$$
where $e_\ga =\int_G \ga (-t)f(t) \, dt$. Thus $e_\ga \cd \ga \in V_f$ for every
$\ga \in \Ga$.
Each $\ga \in \Ga$ is an exponential. Thus linear combinations of the
elements $e_\ga \cd \ga$ are $E$-valued exponential polynomials by Theorem
\ref{t3}. Therefore, it is enough to show that $f$ is in the closure of the set
$N$ of linear combinations of the elements $e_\ga \cd \ga$.

The closure $\cl N$ of $N$ is a closed subspace of $C(G,E)$, so if $f\notin
\cl N$, then there is a $\La \in C(G,E)^*$ such that $\La (f)\ne 0$ and
$\La (e_\ga \cd \ga )=0$ for every $\ga \in \Ga$. Let $M(g)=\La (\hat g )$,
where $\hat g (x)=g(-x)$. Then $(M*g)(0)=\La (g)$ for every $g\in C(G,E)$.
Therefore, we have $M*f \ne 0$ and $M*(e_\ga \cd \ga )=0$ for every
$\ga \in \Ga$. Then
$$(M*f )*\ga =M*(f*\ga )=M*(e_\ga \cd \ga )=0$$
for every $\ga \in \Ga$. Since $M*f\in C(G,\cc )$ and the set of
trigonometric polynomials is dense in $C(G,\cc )$, it follows that
$M*f=0$, a contradiction. \hfill $\square$


\begin{thebibliography}{100}

\bibitem{AD} J. Acz\'el and J. Dhombres: {\it
Functional equations in several variables.} 
Encyclopedia of Mathematics and its Applications, 31. Cambridge University
Press, Cambridge, 1989.

\bibitem{BN} S. Bochner and J. von Neumann, Almost periodic functions in
groups. II, {\it Trans. Amer. Math. Soc.} {\bf 37} (1935), 21-50.

\bibitem{D} J. Diestel: {\it Sequences and Series in Banach Spaces.}
Graduate Texts in Mathematics, 92. Springer-Verlag, New York, Berlin,
Heidelberg, Tokyo, 1984.  

\bibitem{Dj} D. \v Z. Djokovi\'c, A representation theorem for
$(X_1 -1)(X_2 -1) \ldots (X_n -1)$ and its applications,
{\it Ann. Polon. Math.} {\bf 22} (1969), 189-198.

\bibitem{108} M. Laczkovich, Polynomial mappings on Abelian groups, 
{\it Aequationes Math.} {\bf 68} (2004), 177-199.

\bibitem{L} M. Laczkovich, Local spectral synthesis on Abelian groups,
{\it Acta Math. Hungar.} {\bf 142} (2014), 313-329. 

\bibitem{139} M. Laczkovich, A characterization of generalized exponential
polynomials in terms of decomposable functions, {\it Acta Math. Hungar.}
{\bf 158} (2) (2019), 338-351.

\bibitem {140} M. Laczkovich, The Levi-Civita equation in function classes,
to appear in {\it Aequationes  Math.}

\bibitem{LSz1} M. Laczkovich and G. Sz\'ekelyhidi, Harmonic analysis on discrete
Abelian groups, {\it Proc. A.M.S.} {\bf 133} (2005), no. 6, 1581-1586.

\bibitem{LSz2} M. Laczkovich and L. Sz\'ekelyhidi, Spectral synthesis on
discrete Abelian groups, {\it Math. Proc. Camb. Phil. Soc.} {\bf 143} (2007),
103-120.
  
\bibitem{Mc} M. A. McKiernan, Equations of the form $H(x\circ y) =\sum_i f_i
(x)g_i (y),$ {\it Aequationes Math.} {\bf 16} (1977), 51-58.

\bibitem{Ru} W.~Rudin: {\it Functional Analysis.} McGraw-Hill, 1973.
  
\bibitem{Sz1} L. Sz\'ekelyhidi, Note on exponential polynomials,
{\it Pacific J. Math.} {\bf 103} (1982), 583-587.

\bibitem{Sz2} L. Sz\'ekelyhidi: {\it Convolution Type Functional Equations on
Topological Abelian Groups}. World Scientific, 1991.

\bibitem{Sz3}  L. Sz\'ekelyhidi, The failure of spectral synthesis
on some types of discrete Abelian groups. {\it J. Math. Anal. and
Applications} {\bf 291} (2004), 757-763.



\end{thebibliography}
\end{document}